\documentclass[titlepage,11pt]{article}
\usepackage{latexsym}
\usepackage{amsfonts}
\usepackage{amsmath}
\usepackage{amsthm}
\newcommand\blackslug{\hbox{\hskip 1pt \vrule width 4pt height 8pt depth 1.5pt
        \hskip 1pt}}
\newcommand\bbox{\hfill \quad \blackslug \bigbreak}
\def\DD{\hbox{-}}
\def\CC{\hbox{-}\cdots\hbox{-}}
\def\LL{,\ldots,}
\def\cupcup{\cup\cdots\cup}

\title{Induced paths in graphs without anticomplete cycles}
\author{
Tung Nguyen\\
Princeton University, Princeton, NJ 08544
\\
\\
Alex Scott\thanks{Research supported by EPSRC grant EP/X013642/1.}\\
Mathematical Institute, University of Oxford, Oxford OX2 6GG, UK
\\
\\
Paul Seymour\thanks{Supported by AFOSR grant
A9550-19-1-0187.}\\
Princeton University, Princeton, NJ 08544
}

\date{}

\swapnumbers
\newtheorem{thm}{Theorem}[section]

\newtheorem{conj}[thm]{Conjecture}

\newcommand{\Proof}{\noindent{\bf Proof.}\ \ }

\begin{document}
\maketitle
\begin{abstract}
Let us say a graph is {\em $s\mathcal{O}$-free}, where $s\ge 1$ is an integer,  if there do not exist $s$ cycles of the graph that are 
pairwise vertex-disjoint and have no edges joining them.
The structure of such graphs, even when $s=2$, is not well understood. For instance, until now  we did not know how to test 
whether a graph is $2\mathcal{O}$-free in polynomial time; and there was an open
conjecture, due to Ngoc Khang Le, that $2\mathcal{O}$-free graphs have only a polynomial number of induced paths.

In this paper we prove Le's conjecture; indeed, we will show that for all $s\ge 1$, there exists $c>0$ such that 
every $s\mathcal{O}$-free graph $G$ has at most $|G|^c$ induced paths. This provides a poly-time algorithm to test if 
a graph is $s\mathcal{O}$-free, for all fixed $s$.

The proof has three parts. First, there is a short and beautiful proof,
due to Le, that reduces the question to proving the same thing for graphs
with no cycles of length four.
Second, there is a 
recent result of Bonamy, Bonnet, D\'epr\'es, Esperet, Geniet, Hilaire, Thomass\'e and Wesolek, that in every $s\mathcal{O}$-free
graph $G$ with no cycle of length four, 
there is a set of vertices that intersects every cycle, with size logarithmic in $|G|$. And third, there is an argument
that uses the result of Bonamy et al. to deduce the theorem. The last is the main content of this paper.

\end{abstract}

\section{Introduction}
Graphs in this paper are finite and simple (we will occasionally need parallel edges, but then we speak of ``multigraphs'').
Two subsets $X,Y$ of the vertex set of a graph $G$ are {\em anticomplete} if they are disjoint and there is no edge of $G$
between $X$ and $Y$; and we say two subgraphs of $G$ are {\em anticomplete} if their vertex sets are anticomplete.
If $s\ge 1$ is an integer, a graph $G$ is {\em $s\mathcal{O}$-free} if no $s$ cycles of $G$ are pairwise vertex-disjoint and 
anticomplete.
We do not understand such graphs very well: for instance, until now we did not know a polynomial-time algorithm to 
recognize $2\mathcal{O}$-free graphs. In an attempt to find such an algorithm, several years ago Ngoc Khang Le proposed
the (unpublished) conjecture~\cite{khang} that there exists $c>0$ such that every $2\mathcal{O}$-free graph $G$ has at most $|G|^c$ induced cycles; and
the stronger conjecture that the same is true for paths, that is:
\begin{conj}\label{leconj}
There exists $c>0$ such that every $2\mathcal{O}$-free graph $G$ has at most $|G|^c$ induced paths.
\end{conj}
If \ref{leconj} is true, it is easy to derive a poly-time algorithm to test for being $2\mathcal{O}$-free.
Here is a sketch of such an algorithm:
\begin{itemize}
\item For each vertex $v$, find all the induced paths with first vertex $v$. 
\item For each induced path $P$, find all induced cycles that consist of $P$ and one extra vertex.
\item Check whether any two of these cycles are disjoint and have no edges betweeen them.
\end{itemize}
The total running time is proportional to the number of induced paths (times some polynomial in $|G|$), and so is polynomial,
by \ref{leconj}.

In this paper we will prove the stronger conjecture (and hence both conjectures) for $2\mathcal{O}$-free graphs, and indeed 
for $s\mathcal{O}$-free graphs.
More exactly:
\begin{thm}\label{realpolypaths}
For all integers $s\ge 1$ there exists $c>0$ such that, if $G$ is $s\mathcal{O}$-free,
then $G$ has at most $|G|^c$ induced paths.
\end{thm}
(We remark that by definition, every path has at least one vertex; so the one-vertex graph has only one induced path.)
We will use a short and elegant argument, due to Ngoc Khang Le, that deduces \ref{realpolypaths} from the following weaker result:

\begin{thm}\label{polypaths}
For all integers $s\ge 1$ there exists $c>0$ such that, if $G$ is $s\mathcal{O}$-free, and has no cycle of length four,
then $G$ has at most $|G|^c$ induced paths.
\end{thm}

A subset $Z\subseteq V(G)$ is {\em cycle-hitting} if every cycle of $G$ has a vertex in $Z$.
We will show that:
\begin{thm}\label{mainthm}
Let $s\ge 1$ be an integer; then there exist $c_1,c_2,c_3$ such that if $G$ is $s\mathcal{O}$-free, and $Z\subseteq V(G)$ is
a cycle-hitting set, then $G$ has at most
$|G|^{c_1}2^{c_2|Z|+c_3}$ induced paths.
\end{thm}
If $G$ has no cycle of length four then it does not contain $K_{2,2}$ as a subgraph; so to complete the proof of 
\ref{polypaths} and hence of \ref{realpolypaths},
we will use the $t=2$ case of a recent result due to Bonamy, Bonnet, D\'epr\'es,
Esperet, Geniet, Hilaire, Thomass\'e and Wesolek~\cite{frenchpaper}, the following:
\begin{thm}\label{loghittingset}
For all integers $s,t\ge 0$, there exists $c>0$ such that if $G$ is $s\mathcal{O}$-free and does not contain $K_{t,t}$ as a 
subgraph, there is a cycle-hitting set of cardinality at most $c\log|G|$.
\end{thm}

Clearly \ref{polypaths} follows from \ref{mainthm} and \ref{loghittingset}, so the main goal of this paper is to prove \ref{mainthm}.
Let us sketch the idea of its proof.
Let $G$ be an $s\mathcal{O}$-free graph, and let $Z\subseteq V(G)$ be a cycle-hitting set. 
Thus $G\setminus Z$ is a forest $F$ say.
It suffices to count the number of induced paths 
$P$ of $G$ with both ends in $Z$ and with $Z\subseteq V(P)$; because then we can bound the total number of induced paths $P$
by enumerating all possibilities for $V(P)\cap Z$, and for each one, deleting the vertices in $Z\setminus V(P)$, and enumerating 
all possibilities for the two minimal subpaths of $P$ between an end of $P$ and $Z$. So we will focus on such paths $P$,
which we call ``$Z$-covering''. If we want to bound the number of $Z$-covering paths, we can delete any vertices with at least 
three neighbours in $Z$; and we can arrange that $Z$ is stable, by contracting any edges with both ends in $Z$. 
(The number of $Z$-covering paths does not decrease under such contraction, although it might increase.) We need to be careful
with vertices in $V(G)\setminus Z$ that have two neighbours in $Z$, and we will treat such vertices separately. For this sketch, 
let us assume that every vertex in $V(G)\setminus Z$ has at most one neighbour in $Z$, that is, $(G,Z)$ is ``monic''. 
Let $N$ be the set of vertices in $F$ with a neighbour in $Z$.
We are interested in paths of $F$ that join distinct vertices in $N$ and have no internal vertices in $N$ (we call them ``transitions'').
For each transition there are two vertices in $Z$ adjacent to the end of the path (its ``feet''), or maybe only one
such vertex, if it is adjacent to both ends of the path. 
We will show that, by deleting a bounded number of vertices in $F$, and deleting the neighbours of a bounded number of vertices in $Z$,
we can arrange that every surviving transition has two feet, and at most constantly many of them have the same two feet.  
(And it suffices to count the $Z$-covering paths in the part of the graph that survives.)
Next we show (not quite; we will explain later) that we can choose a ``normal'' set of transitions with 
cardinality proportional to $|N|$ (``normal''
means basically that any two of the transitions that are not anticomplete have a common end). But now look at the multigraph
with vertex set $N$ defined by the pairs of feet of the members of the normal set. We can show that this multigraph does not have 
$s$ vertex-disjoint cycles;
because if it does, then $G$ would have $s$ anticomplete cycles (this is why we wanted the set to be normal; this statement is not true
for general sets of transitions, but it works for normal sets). There is a theorem of Erd\H{o}s and P\'osa that says that in such a graph,
there is a set of vertices of bounded size that meets all cycles; so there exists $X\subseteq Z$ of bounded size such that only 
$|Z|$ of the transitions
in the normal set have no foot in $X$. The number that do have a foot in $X$ is also only some constant times $Z$, since only
a bounded number 
have the same pair of feet; so the normal set has cardinality $O(|Z|)$. But its cardinality was proportional to $|N|$, 
and this tells us that $|N|\le O(|Z|)$, and so there are only $ O(|Z|)$ edges between $Z,N$.
Each $Z$-covering path is determined by the set of edges between $Z,N$ that it uses, and there are only $2^{O(|Z|)}$ such subsets, 
so there are only $2^{O(|Z|)}$ $Z$-covering paths, which is what we wanted to show.

Except we cheated in the above; our claim that we can find a large normal set
of transitions is not actually true. What is true is that we can find such a set with size proportional to the number of vertices in $N$
that belong to components of $F$ that have at least two vertices in $N$. We need a special argument to dispose of components of $F$
that only contain one vertex in $N$ that we do not describe here. We also cheated in assuming that $(G,Z)$ is monic, but
the argument we sketched above is the basic idea, and it just needs a few technical patches to make it work.

We recently learned that~\ref{polypaths} was proved independently in the
unpublished paper~\cite{khang}; however, we are told that there are currently no
plans to publish it.

\section{Reducing \ref{realpolypaths} to \ref{polypaths}}

In this section we give the beautiful argument of Ngoc Khang Le, that reduces \ref{realpolypaths} to \ref{polypaths}.
Our thanks to Le for allowing us to include this proof.

Let us say an {\em ordered} induced path is an induced path with one end distinguished as its {\em first} vertex.
We will show the following:
\begin{thm}\label{reduction}
Let $s,c\ge 1$, and suppose that every $s\mathcal{O}$-free graph $G$ with no cycle of length four has at most 
$|G|^{c}$ ordered induced paths.
Then every $s\mathcal{O}$-free graph $G$ has at most $|G|^{d}$ ordered induced paths, where $d=2+ (s-1)(c+6)$.
\end{thm}
\Proof
For $1\le r\le s$, let $d_r= 2+ (r-1)(c+6)$.
We prove a stronger statement, that for $1\le r\le s$, every $r\mathcal{O}$-free graph $G$ has at most   
$|G|^{d_r}$ ordered induced paths.
We proceed by induction on $r$. If $r=1$ then $G$ is a forest, and so has only at most $|G|^2$ ordered induced paths and the 
claim is true. So we assume that $r\ge 2$ and the claim holds for $r-1$. A {\em 4-cycle} means a cycle of length four.

Let $|G|$ be some $r\mathcal{O}$-free graph.
For each 4-cycle $C$ of $G$, let $X_C$ be the set of all vertices of $G$ that are not in $V(C)$ and have no neighbour in $V(C)$.
Let $P$ be an ordered induced path of $G$, with vertices $p_1\CC p_k$ in order, where $p_1$ is its first vertex. 
If there is a 4-cycle $C$
with $p_1\in X_C$, then we may choose $j\in \{1\LL k\}$ maximum such that
there is a 4-cycle $C$ with $p_1\LL p_j\in X_C$, and we call the path $p_1\CC p_j$ the {\em head} of $P$.

Let us first count the number of choices of $P$ that have no head.
There are at most $|G|^2$ choices of $P$ with $k\le 2$, so let us assume $k\ge 3$. There are only $|G|^2$ choices for $p_1$ 
and $p_2$; let us fix some choice of 
$p_1,p_2$, and let $Y$ be the set of vertices different from and nonadjacent to $p_1$. Thus $p_3\CC p_k$ is an ordered induced
path of $G[Y]$; but since $P$ has no head, it follows that $G[Y]$ has no cycle of length four, and so there are only at 
most $|Y|^{c}\le (|G|-2)^{c}$ choices for $p_3\CC p_k$.
Hence altogether there are at most
$|G|^2(|G|-2)^{c}+|G|^2\le |G|^{c+2}$ choices of ordered induced paths $P$ with no head.

Now let us count the number of choices of $P$ that have a head. If some ordered induced path $p_1\CC p_j$ is the head of some 
ordered induced path $P$,
then there is a 4-cycle $C$ such that $p_1\CC p_j$ is a path of $G[X_C]$. But $G[X_C]$ is $(r-1)\mathcal{O}$-free, and 
so, from the inductive hypothesis, it contains at most $|X_C|^{d_{r-1}}\le (|G|-1)^{d_{r-1}}$ ordered induced paths; 
and there are at most
$\binom{|G|}{4}\le (|G|-1)^4$ choices for $C$.
Consequently there are at most 
$(|G|-1)^{d_{r-1}+4}$ choices for the head $p_1\CC p_j$. 

For each choice of head $p_1\CC p_j$, let us count the number of ordered induced paths $p_1\CC p_k$ with head $p_1\CC p_j$.
There are at most $|G|^2$ with $k\le j+2$,
so let us assume that $k\ge j+3$. Again, there are only $|G|^2$
choices for $p_{j+1}$ and $p_{j+2}$; having selected them, let us count the possibilities for $p_{j+3}\CC p_k$.
Let $Z$ be the set of vertices of $G$ different from and nonadjacent to all of $p_1\LL p_{j+1}$. From
the maximality of $j$ in the definition of a head, $G[Z]$ has no 4-cycle, and $p_{j+3}\CC p_k$ is an ordered induced path of $G[Z]$. 
Consequently, having selected $p_1\CC p_{j+2}$, 
there are only $(|G|-2)^{c}$ choices for $p_{j+3}\CC p_k$ with $k\ge j+3$. Hence, having selected $p_1\CC p_j$,
there are at most $|G|^2(1+(|G|-2)^{c})\le |G|^{c+2}$ 
choices for
$p_{j+1}\CC p_k$. Thus, altogether 
there are at most $(|G|-1)^{d_{r-1}+4}|G|^{c+2}$ choices for $p_1\CC p_k$ that have a head.
Including the paths with no head, we have a total of at most
$$|G|^{c+2} + (|G|-1)^{d_{r-1}+4}|G|^{c+2}\le |G|^{c+d_{r-1}+6}=|G|^{d_{r}}$$
ordered induced paths in $G$. This proves \ref{reduction}.~\bbox

\section{Some lemmas about forests}

The remainder of the paper is devoted to proving \ref{mainthm}.
We will need several lemmas about collections of subtrees in a forest.
We begin with
\begin{thm}\label{ring}
Let $F$ be a forest, let $T_1\LL T_{\ell}$ be trees of $F$, and let $H$ be the graph with vertex
set $\{1\LL \ell\}$ in which $i,j$ are adjacent in $H$ if and only if $T_i, T_j$ are not anticomplete.
If $H$ is bipartite then $H$ is a forest.
\end{thm}
\Proof
Since $H$ is bipartite, we may assume that for some $k\in \{0\LL \ell\}$, $T_1\LL T_k$ are pairwise anticomplete, and $T_{i+1}\LL T_\ell$
are pairwise anticomplete.
Suppose that $H$ has a cycle $C$. We may assume that $1,2\in V(C)$.
Let $P_1,P_2$ be the two paths of $C$ between $1,2$. For $h = 1,2$ let $I_h= \{k+1\LL \ell\}\cap V(P_h)$. Let $v_i\in V(T_i)$
for $i = 1,2$. For $h = 1,2$, there is a path $Q_h$ of $F$ between $v_1,v_2$ with interior included in the union of the
sets $V(T_i)\;(i\in V(P_h))$, and hence included in
$$V(T_1\cupcup T_k)\cup \bigcup_{i\in I_h}V(T_i).$$
Since $F$ is a forest, it follows that $Q_1=Q_2$, and so every vertex of $Q_1$ not in $V(T_1\cupcup T_k)$
belongs to both $\bigcup_{i\in I_1}V(T_i)$ and to $\bigcup_{i\in I_2}V(T_i)$, which is impossible since these two sets are disjoint.
Consequently $V(Q_1)\subseteq V(T_1\cupcup T_k)$, which is also impossible since $T_1\LL T_k$ are anticomplete, and $Q_1$
has an end in $T_1$ and an end in $T_2$. This proves \ref{ring}.~\bbox

The next result is related to a result (Theorem 7) of~\cite{kaminski}:
\begin{thm}\label{matroid}
Let $H$ be a forest, let $(A,B)$ be a bipartition of $H$ with $|A|=|B|$, and let i$n$ be an integer with $0\le n\le |A|$. 
Then there is a stable set $X$ of $H$
with $|X|=|A|$ and with $|X\cap A|=n$.
\end{thm}
\Proof
We may assume that $1\le n\le |A|-1$, because otherwise we may take $X\in \{A,B\}$. 
We use induction on $|A|$. Let $v\in V(H)$ have degree at most one. From the symmetry we may assume that
$v\in B$; let $u\in A$ be the neighbour of $v$, if there is one, and otherwise choose $u\in A$ arbitrarily. 
Let $A'=A\setminus \{u\}$, and $B'=B\setminus \{v\}$. 
From the inductive hypothesis, there is a stable set $X'\subseteq A'\cup B'$ with $|X|=|A'|$ and with $|X\cap A'|=n$.
But then $X'\cup \{v\}$ satisfies the theorem. This proves \ref{matroid}.~\bbox

These are used to prove the following:
\begin{thm}\label{selftrees}
Let $F$ be a forest, let $k,s\ge 0$ be integers, and for $1\le i\le s$ let $\mathcal{F}_i$ be a set of $s!k$ paths of $F$, pairwise
anticomplete. Then there exist $P^i_1\LL P^i_k\in \mathcal{F}_i$ for $1\le i\le s$, such that these $sk$ paths are pairwise anticomplete.
\end{thm}
\Proof
We use induction on $s$. For $1\le i\le s$ let $A_i$ be the set of all pairs $(i,j)$ with $1\le j\le s!k$; and let
$H$ be the graph with vertex set $A_1\cupcup A_s$, where
$(i,j)$ and $(i',j')$ are adjacent if the $j$th member of $\mathcal{F}_i$ is not anticomplete to the $j'$th member of $\mathcal{F}_{i'}$.
By \ref{ring} applied to $\mathcal{F}_1$ and $\mathcal{F}_2$, for $2\le j\le s$ the subgraph of $H$
induced on $A_1\cup A_i$ is a forest, with a bipartition
$(A_1,A_i)$; and by \ref{matroid},
there is a stable set $X_{i}$ of $H$ with cardinality
$s!k$, containing
$(s!-(s-1)!)k$ vertices of $A_1$ and $(s-1)!k$ vertices of $A_i$.  The sets $A_1,X_{2}\LL X_{s}$ have at least $(s-1)!k$ vertices 
in common; and so for $1\le i\le s$ there is a subset $\mathcal{F}'_i$ of $\mathcal{F}_i$
with cardinality $(s-1)!k$, such that all the paths in $\mathcal{F}'_1$ are anticomplete to all the paths in
$\mathcal{F}'_i$ for $2\le i\le s$. But then the result follows from the inductive hypothesis applied to the sets
$\mathcal{F}'_i$ for $2\le i\le s$. This proves \ref{selftrees}.~\bbox

We will also need:

\begin{thm}\label{subtrees}
Let $F$ be a forest, let $n\ge 0$ be an integer, and let $T_1\LL T_k$ be trees of $F$. 
\begin{itemize}
\item If no $n$ of $T_1\LL T_k$ are pairwise vertex-disjoint, there exists $X\subseteq V(F)$ with $|X|\le n-1$ such that
$X\cap V(T_i)\ne \emptyset$ for $1\le i\le k$;
\item If no $n$ of $T_1\LL T_k$ are pairwise anticomplete, there exists $X\subseteq V(F)$ with $|X|\le 2(n-1)$ such that 
$X\cap V(T_i)\ne \emptyset$ for $1\le i\le k$.
\end{itemize}
\end{thm}
\Proof
The first claim is well-known and easy, and we assume it without proof. For the second, let $F'$ be the forest obtained from
$F$ by subdividing once each edge $e$ of $F$ (let $v_e$ be the new vertex that subdivides $e$). For $1\le i\le k$, 
let $T_i'$ be the tree of $F'$ induced on the union of $V(T_i)$ and the set of all $v_e$ such that $e\in E(F)$ has an end in $V(T_i)$.
The hypothesis implies that no $n$ of $T_1'\LL T_k'$ are pairwise vertex-disjoint, and so the result follows by applying the 
first bullet of the theorem to $F'$ and $T_1'\LL T_k'$. This proves \ref{subtrees}.~\bbox

\section{Plantations and transitions}

Let $G$ be an $s\mathcal{O}$-free graph, and let $Z\subseteq V(G)$ be a cycle-hitting set.
We call $(G,Z)$ a {\em plantation}.  (So the definition of a plantation depends on $s$, but we leave this implicit: $s$ will be fixed
throughout anyway.)
Let $F$ be the forest $G\setminus Z$, and let $N$ be the set of vertices in $V(G)\setminus Z$ with a neighbour in $Z$.
We say $(G,Z)$ is  {\em monic} if $Z$ is stable and each vertex in $N$ has a unique neighbour in $Z$.
Let us say a {\em transition} of $(G,Z)$ is a path of $F$ of length at least one, with both ends in $N$ and with no internal vertex in $N$.
Let $P$ be a transition. If $z\in Z$ is adjacent to an end of $P$, we say $z$ is a {\em foot} of $P$. If $(G,Z)$ is monic,
every transition $P$ has one or two feet, and these are the only vertices in $Z$ that have a neighbour in $V(P)$.
We remark that distinct transitions cannot have the same pair of ends, since $F$ is a forest, but they may have the same pair of feet.
If $P$ only has one foot, $P$ is a {\em self-transition}.
We say $(G,Z)$ is {\em selfless} if there is no self-transition.
Starting with a monic plantation, our first objective is to eliminate self-transitions.

We will use two operations to eliminate self-transitions: deletion and explosion. If $(G,Z)$ is a plantation, and $v\in V(G)\setminus Z$, 
then $(G\setminus \{v\},Z)$ is a plantation, monic if $(G,Z)$ is monic. 
Moreover, each transition of $(G\setminus \{v\},Z)$ 
is a transition of $(G,Z)$,
so deleting vertices
in $V(G)\setminus Z$ may be used to eliminate some self-transitions, without introducing new ones. 
Second, if $v\in Z$, let $G'$ be obtained from $G$ by deleting $v$ and all its neighbours in $V(G)\setminus Z$.
Then again $(G',Z\setminus \{z\})$ is a plantation, monic if $(G,Z)$ is monic, and each of its transitions is a transition 
of $(G,Z)$. This operation is called
{\em exploding} $v$. We will show:

\begin{thm}\label{selfless}
Let $(G,Z)$ be a monic plantation. Then there exist $X\subseteq Z$ and 
$Y\subseteq V(G)\setminus Z$, with $|X|<s$ and $|Y|<2s\cdot s!$, such that exploding the vertices in $X$
and deleting the vertices in $Y$ yields a selfless plantation.
\end{thm}
\Proof
As before, let $F=G\setminus Z$, and let $N$ be the set of vertices in $V(G)\setminus Z$ with a neighbour in $Z$.
Let us say $z\in Z$ is {\em $k$-self-important} if there are $k$ self-transitions, pairwise anticomplete and each with foot $z$. 
\\
\\
(1) {\em There do not exist $s$ distinct vertices in $Z$ that are $s!$-self-important.}
\\
\\
Suppose that $z_1\LL z_s\in Z$ are each $s!$-self-important, and for $1\le i\le s$ let $\mathcal{F}_i$
be a set of $s!$ self-transitions, each with foot $z$ and pairwise anticomplete.
By \ref{selftrees} with $k=1$, there exist $P_i\in \mathcal{F}_i$ for $1\le i\le k$, such that $P_1\LL P_s$
are pairwise anticomplete. Thus $V(P_i)\cup \{z_i\}$ induces a cycle $C_i$ say, for each $i$, and since $(G,Z)$ is monic, 
$z_i$ has no neighbour in $C_j$ if $i,j$ are distinct, and so $C_1\LL C_s$ are pairwise anticomplete, a contradiction.
This proves (1).
\\
\\
(2) {\em If there is no $s!$-self-important vertex in $Z$,
then there exists $Y\subseteq V(F)$ with $|Y|\le 2s\cdot s!$ such that 
deleting the vertices in $Y$ yields a selfless plantation.}
\\
\\
We claim that there do not exist $s\cdot s!$ self-transitions that are pairwise anticomplete; for if there are, then since 
no $s!$ of them have the same foot, we could choose $s$ of them all with distinct feet (each with only one foot, but all distinct);
and again that gives us $s$ pairwise anticomplete cycles, a contradiction.
From \ref{subtrees}, there exists $Y\subseteq V(F)$ with $|Y|<2s\cdot s!$ 
such that every self-transition contains a vertex in $Y$; and so deleting the vertices in $Y$ yields a selfless plantation.
This proves (2).

\bigskip
But from (1), by exploding at most $s-1$ vertices in $Z$, we can produce a plantation with no $s!$-self-important vertex; and so 
the result follows from (2). This proves \ref{selfless}.~\bbox

If $P$ is a path, we denote the interior of $P$ (that is, the set of vertices that have degree two in $P$) by $P^*$.
Let $(G,Z)$ be a monic selfless plantation. 
If $z,z'\in Z$, the {\em multiplicity} of the pair $(z,z')$ is the number of transitions with feet $z,z'$.
Thus the multiplicity of $(z,z)$ is zero, since $(G,Z)$ is selfless.
We say that $(G,Z)$ has {\em thickness} $k$ if $k$ is the maximum of the multiplicity of pairs of elements of $Z$.
Our next objective is to obtain a plantation with bounded thickness, again by deleting and exploding a bounded number of vertices. 
We will show the following.
\begin{thm}\label{getthin}
Let $(G,Z)$ be a monic selfless plantation. Then there exists $X\subseteq Z$ with
$|X|\le 6s-4$ such that
exploding the vertices in $X$
yields a plantation with thickness at most $2\cdot s!(2\cdot s!+s)$.
\end{thm}
\Proof
Let $N$ be the set of vertices in $V(G)\setminus Z$ with a neighbour in $Z$, and let $F$ be the forest $G\setminus Z$.
We observe first:
\\
\\
(1) {\em Let $z,z'\in Z$. If $P_1,P_2$ are distinct transitions both with feet $z,z'$, then $P_1^*, P_2^*$ are anticomplete, and either 
\begin{itemize}
\item $P_1,P_2$ are anticomplete; or
\item $P_1,P_2$ have a common end and $P_1\cup P_2$ is an induced path; or
\item $V(P_1),V(P_2)$ are disjoint and there is a unique edge between them, joining an end of $P_1$ and an end of $P_2$.
\end{itemize}
}
\noindent Let $P_i$ have ends $a_i,b_i$ for $ i = 1,2$, where $a_1,a_2$ are adjacent to $z$, and $b_1,b_2$ to $z'$.
Since $P_1,P_2$ are distinct, and they are both paths in the forest $F$, they do not have the same pairs of ends; and so we may assume
that $a_1\ne a_2$. 
Let $T_1$ be the maximal tree of $F$ that contains $P_1$ and has the property that every vertex in $N\cap V(T_1)$ has degree one in $T_1$.
Since $(G,Z)$ is selfless, $a_1$ is the only neighbour of $z$ in $V(T_1)$,
and so $a_2\notin V(T_1)$; and consequently $P_2^*\cap V(T_1)=\emptyset$. Similarly, either $b_2\notin V(T_1)$ or $b_2=b_1$.
The vertices of $P_1^*$
are not leaves of $T_1$, and so every vertex of $G$ with a neighbour in $P_1^*$ belongs to $V(T_1)$. Consequently
$P_1^*, P_2^*$ are anticomplete, and $a_2$ has no neighbour in $P_1^*$, and $b_2$  has no neighbour in $P_1^*$ unless $b_1=b_2$.
Similarly $a_1$ has no neighbour in $P_2^*$, and $b_1$  has no neighbour in $P_2^*$ unless $b_1=b_2$.

If $V(P_1)\cap V(P_2)\ne \emptyset$, then $P_1,P_2$ have a common end, and so $b_1=b_2$; but then the second outcome holds.
Thus we may assume that $V(P_1), V(P_2)$ are disjoint. If they are anticomplete, then the first outcome holds; and if not, 
the edge between $V(P_1),V(P_2)$ is unique (since $F$ is a forest) and the third outcome holds. This proves (1).

\bigskip

Let $z,z'\in Z$. If $P_1\LL P_k$ are transitions that are pairwise anticomplete, and all with the same feet $z,z'$, we call $\{P_1\LL P_k\}$ a 
{\em $(z,z')$-linkage}. If $P_1\LL P_k$ all have a common end, we call $\{P_1\LL P_k\}$ a
{\em $(z,z')$-star}, and the common end of $P_1\LL P_k$ is called the {\em centre}. 
\\
\\
(2) {\em Let $z,z'\in Z$, let $p,q\ge 0$ be integers, and let $(z,z')$ have multiplicity at least $2pq$. Then 
there is either a  $(z,z')$-linkage of cardinality $p$, or a $(z,z')$-star of cardinality $q$.}
\\
\\
Let $P_i\;(i\in I)$ all be distinct transitions, with the same feet $z,z'$, where $|I|= 2pq$. For each $i\in I$ let $P_i$ have ends $a_i, b_i$,
where $a_i$ is adjacent to $z$ and $b_i$ to $z'$. Every bipartite graph with $2pq$ edges has a matching of size $2p$ or a vertex 
of degree at least $q$, from K\"onig's theorem; and because of this, applied to the bipartite graph with bipartition 
$(\{a_i:i\in I\},\{b_i:i\in I\})$ and edge set $\{\{a_i,b_i\}:i\in I\}$, we may assume that either $a_1\LL a_{2p}, b_1\LL b_{2p}$
are all distinct, or $a_1=\cdots =a_q$. In the second case, $\{P_1\LL P_q\}$ is a
$(z,z')$-star by (1), so we assume the first holds. Let $H$ be the graph with vertex set $\{1\LL 2p\}$, in which $i,j$ are adjacent if 
$P_i,P_j$ are not anticomplete (and hence they are vertex-disjoint and there is a unique edge between them, by (1)). 
A graph isomorphic to $H$
can be obtained from $F$ by deleting all vertices not in $P_1\LL P_{2p}$ and contracting the edges of $P_1\LL P_{2p}$; and so 
$H$ is a forest. Hence it has a stable set of cardinality $p$, say $\{1\LL p\}$; and then 
$\{P_1\LL P_p\}$ is a $(z,z')$-linkage. This proves (2).
\\
\\
(3) {\em There do not exist distinct $z_1,z_1',z_2,z_2'\LL z_s, z_s'\in Z$ 
such that for $1\le i\le s$ there is a $(z_i,z_i')$-linkage of cardinality $2 \cdot s!$.}
\\
\\
Suppose such vertices exist, and for $1\le i\le s$ let $\mathcal{F}_i$ be a set of $2\cdot s!$ transitions each with feet $z_i, z_i'$, 
and pairwise
anticomplete. By \ref{selftrees} with $k=2$, for $1\le i\le s$ there exist distinct $P_i,Q_i\in \mathcal{F}_i$
such that $P_1,Q_1\LL P_s,Q_s$ are pairwise anticomplete. But then the cycles induced on $V(P_i)\cup V(Q_i)\cup \{z_i,z_i'\}$
are pairwise anticomplete, a contradiction. This proves (3).
\\
\\
(4) {\em  There do not exist distinct $z_1,z_1',z_2,z_2'\LL z_{2s}, z_{2s}'\in Z$
such that for $1\le i\le {2s}$ there is a $(z_i,z_i')$-star of cardinality $2\cdot s!+s$.}
\\
\\
Suppose such vertices exist. The centres of the $2s$ stars are distinct vertices of $F$, and hence some $s$ of them are pairwise nonadjacent;
thus we may assume that $\mathcal{S}_i$ is a $(z_i, z_i')$-star of cardinality $2\cdot s!+s$ with centre $a_i$ for $1\le i\le s$, and $a_1\LL a_s$ are pairwise 
nonadjacent. Let $i,j\in \{1\LL s\}$ be distinct. Since $a_i\in N$, it does not belong to the interior of any member of $\mathcal{S}_j$;
and since $z_1,z_1',z_2,z_2'\LL z_{s}, z_{s}'\in Z$ are distinct and $(G,Z)$ is monic, $a_i$ is not an end of any member of $\mathcal{S}_j$.
Since $F$ is a forest, $a_i$ has a neighbour in at most one member of $\mathcal{S}_j$. Thus for $1\le i\le s$, there are at most $s-1$
members of $\mathcal{S}_i$ that contain a neighbour of $a_j$ for some $j\in \{1\LL s\}\setminus \{i\}$; and so we may choose
$\mathcal{S}'_i\subseteq \mathcal{S}_i$ of cardinality $2\cdot s!$ such that no member of $\mathcal{S}'_i$ contains any vertex adjacent to
some $a_j$ with $j\ne i$. For each $P\in \mathcal{S}_i$, let us say $P\setminus \{a_i\}$ is its {\em truncation}; and let $\mathcal{F}_i$
be the set of truncations of the members of $\mathcal{S}_i'$. Thus the members of $\mathcal{F}_i$ are pairwise anticomplete.
By \ref{selftrees} with $k=2$, there exist distinct $Q_i, Q_i'\in \mathcal{F}_i$ for $1\le i\le s$, such that $Q_1,Q_1'\LL Q_s, Q_s'$
are pairwise anticomplete. But for $1\le i\le s$, there is a cycle $C_i$ with $V(C_i)\subseteq V(Q_i)\cup V(Q_i')\cup \{a_i, z_i, z_i'\}$,
and these $s$ cycles are pairwise anticomplete, a contradiction. This proves (4).

\bigskip

Choose distinct $z_1,z_1',z_2,z_2'\LL z_{r}, z_{r}'\in Z$ with $r$ maximum such that for $1\le i\le r$ there is a $(z_i,z_i')$-linkage 
of cardinality $2\cdot s!$. Let $X_1=\{z_1,z_1',z_2,z_2'\LL z_{r}, z_{r}'\}$. From (3), $r\le s-1$, and so $|X_1|\le 2(s-1)$;
and from the maximality of $r$, 
for all $z,z'\in Z$,
if there is a $(z,z')$-linkage of  
cardinality $2\cdot s!$ then one of $z,z'\in X_1$. Similarly from (4), there is a set $X_2\subseteq Z$ with $|X_2|\le 2(2s-1)$
such that for all $z,z'\in Z$,
if there is a $(z,z')$-star of          
cardinality $2\cdot s!+s$ then one of $z,z'\in X_2$. Hence from (1), for all $z,z'\in Z$,
if $(z,z')$ has multiplicity at least $(2s\cdot !)(2\cdot s!+s)$,          
then one of $z,z'\in X_1\cup X_2$. Thus the plantation produced by exploding the vertices in $X_1\cup X_2$ has thickness at most
$(2\cdot s!)(2\cdot s!+s)$. 
This proves 
\ref{getthin}.\bbox

\section{Applying the Erd\H{o}s-P\'osa theorem}
Let $(G,Z)$ be a plantation; we say a set $\mathcal{S}$ of transitions in $(G,Z)$ is {\em normal} if 
\begin{itemize}
\item for all $P,Q\in \mathcal{S}$,
either $P,Q$ are anticomplete or $P,Q$ have a common end; and
\item for each $P\in \mathcal{S}$, there is an edge $e$ of $P$ that does not belong to any other member of $\mathcal{S}$.
\end{itemize}
We need first:
\begin{thm}\label{gettrans}
Let $(G,Z)$ be a plantation, and let $N$ be the set of vertices in $V(G)\setminus Z$ with a neighbour in $Z$. Suppose that
every component of $F$ contains at least two vertices of $N$. Then there
is a normal set $\mathcal{S}$ of transitions with $|\mathcal{S}|\ge |N|/4$.
\end{thm}
\Proof
Let $F$ be the forest $G\setminus Z$. By choosing transitions from each component of $F$ separately, we may assume that $F$
is a tree, and $|N|\ge 2$. If $|N|\le 3$ the result is clear, so we may assume that $|N|\ge 4$.
Choose some vertex $r\in N$, call it the {\em root} of $F$,
and direct every edge of $F$ towards $r$. Let $\mathcal{R}$ be the set of all transitions of $(G,Z)$ that are directed paths.
Thus $|\mathcal{R}|= |N|-1$, since every vertex in $N$ different from $r$ is the first vertex of a unique directed transition.
Moreover, for the same reason, every member of $\mathcal{R}$ has an edge that does not belong to any other member of $\mathcal{R}$.
We will show that there is a normal subset of $\mathcal{R}$ with cardinality at least $|\mathcal{R}|/3$.

Let $P$ be a directed transition, and let $Q$ be the directed path of $F$ from the first vertex of $P$
to the root of $F$. It follows that $P$ is an initial subpath of $Q$. We define the {\em height} of $P$
to be the number of vertices of $Q$ that belong to $N$.
\\
\\
(1) {\em Let $P_1,P_2$ be directed transitions, with heights $h_1,h_2$ where $h_1-h_2$ is a multiple of three. Then either $P,Q$ are anticomplete, or they 
have the same last vertex and therefore the same height.}
\\
\\
Let $P_i$ have first vertex $a_i$ and last vertex $b_i$ for $i = 1,2$. We may assume that $b_1\ne b_2$, and so 
$V(P_1)\cap V(P_2)=\emptyset$. Hence we may assume that there is an edge of $F$ with one end in $V(P_1)$ and the other in $V(P_2)$,
and we may assume this edge is directed from its end $c_1\in V(P_1)$ to its end $c_2\in V(P_2)$, by exchanging $P_1,P_2$
if necessary. Since $c_1$ has at most one out-neighbour in $F$, and $c_2\notin V(P_1)$, it follows that $c_1=b_1$.
For $i = 1,2$, let $Q_i$ be the directed path of $F$ from $a_i$ to the root of $F$. It follows that the edge $c_1c_2$ belongs 
to $Q_1$, and so $Q_1$ contains all the vertices of $N\cap V(Q_2)$ except possibly $a_2$, and in addition contains $a_1,b_1$.
Thus $h_1-h_2\in \{1,2\}$, contradicting that $h_1-h_2$ is a multiple of three. This proves (1).

\bigskip
For $i=1,2,3$, let $\mathcal{S}_i$ be the set of all directed transitions with height congruent to $i$ modulo three.
By (1), each of these sets is normal, and every directed transition belongs to one of them, so one of them has cardinality
at least $|\mathcal{R}|/3=(|N|-1)/3$, and hence at least $|N|/4$, since $|N|\ge 4$. This proves \ref{gettrans}.~\bbox

We need the following result, a theorem of Erd\H{o}s and P\'osa~\cite{EP}:
\begin{thm}\label{EPthm}
If $s\ge 0$ is an integer, there exists $\phi(s)\ge 0$ with the following property. If $G$ is a multigraph in which no $s$ cycles are pairwise vertex-disjoint, there is a subset 
$X\subseteq V(G)$ with $|X|\le \phi(s)$ such that
every cycle of $G$ contains a vertex in $X$.
\end{thm}
Erd\H{o}s and P\'osa showed there exist $c_1,c_2$ such that $c_1 s \log s\le \phi(s) \le c_2 s \log s$ for all $s$, but that does not matter for us. Through the rest of the paper, we use the notation $\phi(s)$ with its meaning in \ref{EPthm}.

We need anticomplete cycles, not just disjoint cycles: but by selecting some transitions carefully, we can make a derived graph,
disjoint cycles in which would yield anticomplete cycles in the original graph.
We use \ref{EPthm} to show the following:
\begin{thm}\label{applyEP}
Let $(G,Z)$ be a monic plantation,  and let $N$ be the set of vertices in $V(G)\setminus Z$
with a neighbour in $Z$. Let $\mathcal{S}$ be a normal set of transitions. Then there exists $X\subseteq Z$ with $|X|\le \phi(s)$
such that at most $|Z|$ members of $\mathcal{S}$ have no foot in $X$.
\end{thm}
\Proof
Let $H$ be the multigraph with vertex set $Z$, edge set $\mathcal{S}$, and incidence relation defined as follows:
for each $P\in \mathcal{S}$, and each $z\in Z$, $P$ is incident with $z$ in $H$ if $z$ is a foot of $P$.
We observe:
\\
\\
(1) {\em If $C$ is a cycle of $H$, there is a cycle $C'$ of $G$ with $V(C')\cap Z\subseteq V(C)$, and $V(C')\setminus Z$
is a subset of the union of the vertex sets of the transitions in $E(C)$.}
\\
\\
Let the vertices and edges of $C$ in order be
$u_1,P_1,u_2,P_2\LL u_m, P_m, u_{m+1}=u_1$. 
Thus $u_1\LL u_m\in Z$ are distinct, and for $1\le i\le m$, $P_i\in \mathcal{S}$ is a
transition with feet $u_i, u_{i+1}$, and $P_1\LL P_m$ are all distinct. 
Suppose that $m=1$; then $H$ has a loop $P_1$, incident with $u_1$ in $H$. Let $p,q$ be the ends of the path $P_1$ in $G$; then the union of $P_1$ with the path $p\DD u_1\DD q$ is the desired cycle. Thus we may assume that $m\ge 2$.

For $1\le i\le m$, let $P_i^+$ be the path between $u_i, u_{i+1}$ with interior $V(P_i)$.
Since $\mathcal{S}$ is normal, there is an edge $e$ of $P_1$ that belongs to none
of $P_2\LL P_m$. But the union of $P_1^+\setminus \{e\}$ and $P_2^+\cupcup P_k^+$ is a connected graph, containing both ends of $e$;
and so contains a path joining the ends of $e$. Adding $e$ to this path gives the desired cycle $C'$. This proves (1).
\\
\\
(2) {\em No $s$ cycles of $H$ are vertex-disjoint.}
\\
\\
Suppose that $C_1\LL C_s$ are $s$ cycles of $H$ that are vertex-disjoint. By (1), there is a
cycle $C_i'$ of $G$ with $V(C_i')\cap Z\subseteq V(C_i)$, and $V(C_i')\setminus Z$
is a subset of the union of the vertex sets of the transitions in $E(C_i)$. 
Since $G$ is $s\mathcal{O}$-free, we may assume that $C_1'$ is not anticomplete to $C_2'$.
Since $C_1, C_2$ are vertex-disjoint, and $Z$ is stable, it follows that $V(C_1')\cap Z$ is anticomplete to $V(C_2')\cap Z$.
Let the vertices and edges of $C_1$ in order be
$$u_1,P_1,u_2,P_2\LL u_m, P_m, u_{m+1}=u_1,$$ 
and define $v_1,Q_1,v_2,Q_2\LL v_n, Q_n, u_{n+1}=v_1$ similarly for $C_2$.
For $1\le i\le m$, two vertices in $V(C_1')\cap Z$ are adjacent to ends of $P_i$, and since $(G,Z)$ is monic, no other
vertices in $Z$ have neighbours in $V(P_i)$. Consequently $V(C_2')\cap Z$ is anticomplete to $V(C_1')$ and similarly
$V(C_1')\cap Z$ is anticomplete to $V(C_2')$. 
Therefore we may assume that $P_1$ is not anticomplete to $Q_1$. Since $\mathcal{S}$ is normal, it follows that
$P_1,Q_1$ have a common end $a$ say; but then the unique neighbour $z\in Z$ of $a$ belongs to both $V(C_1), V(C_2)$,
a contradiction. This proves (2).

\bigskip

From \ref{EPthm}, there exists $X\subseteq Z$ with $|X|\le \phi(s)$ such that $H\setminus X$ is a forest, and therefore has at most
$|Z\setminus X|-1\le |Z|$ edges; and so at most $|Z|$ members of $\mathcal{S}$ have no neighbour in $X$.
This proves \ref{applyEP}.~\bbox

We use this to show:
\begin{thm}\label{Nbound}
Let $(G,Z)$ be a monic selfless plantation, with thickness $k$,  and let $N$ be the set of vertices in $V(G)\setminus Z$
with a neighbour in $Z$. Suppose that every component of $G\setminus Z$ contains at least two vertices in $N$.
Then $|N|\le 4(k\phi(s)+1)|Z|$.
\end{thm}
\Proof 
By \ref{gettrans}, there is a normal set $\mathcal{S}$ of transitions with $|\mathcal{S}|\ge |N|/4$.
From \ref{applyEP}, there exists $X\subseteq Z$ with $|X|\le \phi(s)$ such that 
at most $|Z|$ members of $\mathcal{S}$ have no neighbour in $X$. But since $(G,Z)$
has thickness $k$, for each $x\in X$ and $z\in Z$, there are at most $k$ transitions with feet $x,z$, and therefore for each $x\in X$,
at most $k|Z|$ transitions in $\mathcal{S}$ contain a neighbour of $x$. Since $|X|\le \phi(s)$, it follows that
$|\mathcal{S}|\le k\phi(s)|Z|+|Z|$. But $|\mathcal{S}|\ge |N|/4$, and so $|N|\le 4(k\phi(s)+1)|Z|$.
This proves \ref{Nbound}.~\bbox

\section{Non-monic plantations}
The result \ref{Nbound} brings us close to what we want, but only for monic plantations. In this section we extend it to more
general plantations. Let us say a plantation $(G,Z)$ is {\em dyadic} if $Z$ is stable and 
every vertex in $V(G)\setminus Z$ has at most two neighbours in $Z$.
We say $v\in V(G)\setminus Z$ is {\em binary}
if it has two neighbours in $Z$.
\begin{thm}\label{simpledouble}
Let $(G,Z)$ be a dyadic plantation. Then there exists $X\subseteq Z$ with $|X|\le 2\phi(s)$ such that exploding $X$ yields a dyadic
plantation
with  at most $2|Z|$ binary vertices.
\end{thm}
\Proof
We claim first:
\\
\\
(1) {\em Let $Y$ be a stable set of binary vertices. Then there exists $X\subseteq Z$ with $|X|\le 2\phi(s)$
such that at most $|Z|$ vertices in $Y$ have no neighbour in $X$.}
\\
\\
Let $H$ be the multigraph with vertex set $Z$ and edge set $Y$, where $y\in Y$ is incident in $H$ with $z\in Z$ if $y$ is adjacent
to $z$ in $G$. For every cycle $C$ of $H$, there is a cycle $C'$ of $G$ induced on the vertices of $C$ that are vertices or edges of $C$; and if $C,D$ are vertex-disjoint cycles of $H$, the corresponding cycles $C',D'$ of $G$ are anticomplete (since $Y$ is stable,
$Z$ is stable, and each vertex in $Y$ has exactly two neighbours in $Z$). Consequently no $s$ cycles of $H$ are pairwise vertex-disjoint,
and so
by \ref{EPthm}, there exists $X\subseteq Z$ with $|X|\le \phi(s)$ such that $H\setminus X$ is a forest, and so has at most $|Z|$ edges.
Hence at most $|Z|$
vertices in $Y$ have no neighbour in $X$. This proves (1).

\bigskip

Let $N_2$ be the set of all binary vertices. Since $G\setminus Z$ is a forest and hence bipartite, it follows that $N_2$
is the union of two stable sets; and so by (1) applied to each of these sets, we deduce that there exists $X\subseteq Z$
with $|X|\le 2\phi(s)$ such that at most $2|Z|$ vertices in $N_2$ have no neighbour in $X$. But then $X$ satisfies the theorem.
This proves \ref{simpledouble}.~\bbox

For $z\in Z$, $N(z)$ denotes the set of neighbours of $z$, and for $Z'\subseteq Z$, $N(Z')$ denotes the union of the 
sets $N(z) (z\in Z')$.
We deduce:
\begin{thm}\label{countedges}
Let $(G,Z)$ be a dyadic plantation. Then there exist $X\subseteq Z$ with $|X|\le 2\phi(s)+7s-4$ and $Y\subseteq V(G)\setminus Z$
with $|Y|\le 2s\cdot s!$ and with the following property. Let $F=G\setminus Z$.
For $i = 1,2$, let $N_i$ be the set of all $v\in V(F)\setminus (Y\cup  N(X))$ that have exactly $i$ neighbours in $Z$; 
and let $N_0$ be the set of all $v\in N_1$ such that the component of $F\setminus (Y\cup N(X)\cup N_2)$
containing $v$ contains no other vertex in $N_1$.
Then there are at most 
$$8(s!(2\cdot s!+s)\phi(s)+1)|Z|+ 4s\cdot s!$$ 
edges between $Z\setminus X$ and $V(F)\setminus (N(X)\cup N_0)$.
\end{thm}
\Proof 
By \ref{simpledouble}, there exists $X_1\subseteq Z$ with $|X_1|\le 2\phi(s)$ such that 
exploding $X_1$ yields a dyadic
plantation $(G_1,Z\setminus X_1)$
with  at most $2|Z|$ binary vertices. Let $Y_1$ be the set of binary vertices of $(G_1,Z\setminus X_1)$. It follows that 
$(G_1\setminus Y_1, Z\setminus X_1)$ is monic and $|Y_1|\le 2|Z|$. By \ref{selfless} applied to $(G_1\setminus Y_1, Z\setminus X_1)$, 
there exists 
$X_2\subseteq Z\setminus X_1$ and
$Y\subseteq V(G_1)\setminus (Y_1\cup Z)$, with $|X_2|\le s$ and $|Y|\le 2s\cdot s!$, such that 
starting with $(G_1\setminus Y_1, Z\setminus X_1)$, and exploding the vertices in $X_2$
and deleting the vertices in $Y$, yields a selfless plantation $(G_2,Z\setminus (X_1\cup X_2))$ say.
By \ref{getthin},
there exists $X_3\subseteq Z\setminus (X_1\cup X_2)$ with
$|X_3|\le  6s-4 $ such that starting with $(G_2,Z\setminus (X_1\cup X_2))$ and 
exploding the vertices in $X_3$
yields a monic selfless plantation $(G_3,Z\setminus (X_1\cup X_2\cup X_3))$ with thickness at most $(2\cdot s!)(2\cdot s!+s)$. Let $Y_3$ be the union of the vertex
sets of all components of $G_3\setminus Z$ that have at most one vertex with a neighbour in $Z\setminus (X_1\cup X_2\cup X_3)$.
The plantation $(G_3\setminus Y_3,Z\setminus (X_1\cup X_2\cup X_3))$ satisfies the hypothesis of \ref{Nbound}, and its thickness is at most $(2\cdot s!)(2\cdot s!+s)$, and so
by \ref{Nbound}, there are at most $4(2\cdot s!(2\cdot s!+s)\phi(s)+1)|Z|$
edges between $Z\setminus (X_1\cup X_2\cup X_3)$ and $V(G)\setminus (Y_3\cup Z)$.

Let $X=X_1\cup X_2\cup X_3$; we will show that $X,Y$ satisfy the theorem.  Certainly
$$|X|=|X_1|+|X_2|+|X_3|\le 2\phi(s)+ s+ 6s-4=2\phi(s)+7s-4,$$
and $|Y|\le  2s\cdot s!$.
We recall that $(G_3,Z\setminus X)$ is obtained from $(G,Z)$ by exploding the vertices in $X$ and deleting the vertices
in $Y_1\cup Y$. Let $(G',Z\setminus X)$ be obtained from $(G,Z)$ by exploding the vertices in $X$ and deleting
the vertices in $Y$. 
There are only $2|Y_1|\le 4|Z|$ edges of $G'$ between $Y_1$ and $|Z|$
since $|Y_1|\le 2|Z|$ and each of its members has only two neighbours in $Z$. Thus there are at most
 $4(2\cdot s!(2\cdot s!+s)\phi(s)+2)|Z|$ edges of $G'$ between  $V(G')\setminus (Y_3\cup Z)$ and $Z\setminus X$, that is, 
 between $V(F)\setminus (N(X)\cup Y\cup N_0)$ and $Z\setminus X$. Since $|Y|\le 2s\cdot s!$,
there are only $4s\cdot s!$ edges between $Y$ and $Z$. This proves \ref{countedges}.~\bbox

\section{Counting paths}
Let $(G,Z)$ be a plantation. We denote by $n(G,Z)$ the number of induced paths $P$ of $G$ with $Z\subseteq V(P)$
such that both ends of $P$ belong to $Z$. Let us call such a path $P$ a {\em $Z$-covering path}. Our objective is to show that $n(G,Z)$ is at most the product of a polynomial in
$|G|$ and an exponential in $|Z|$. 

It is enough to work with dyadic plantations, because of the following.
\begin{thm}\label{simplify}
Let $(G,Z)$ be a plantation. Then there is a dyadic plantation $(G',Z')$ with $|G'|\le |G|$ and $|Z'|\le |Z|$
such that $n(G,Z)\le n(G',Z')$.
\end{thm}
\Proof
We prove this by induction on $|G|$. We observe first:
\begin{itemize}
\item If some vertex $v\in V(G)\setminus Z$ has more than two neighbours in $Z$, this vertex does not belong to any $Z$-covering path,
and so we may delete it without changing the number of $Z$-covering paths. Hence in this case we can win by induction on $|G|$; so
we may assume there is no such vertex.
\item If some vertex $v\in V(G)\setminus Z$ has two neighbours $z,z'\in Z$, and $z,z'$ are adjacent, then again $v$ does not belong
to any $Z$-covering path, and we can delete it and win as before. So we may assume that there is no such vertex.
\item If some three vertices in $Z$ are pairwise adjacent, then $n(G,Z)=0$, so we may assume there is no such triangle.
\end{itemize}

If some two vertices $z,z'\in Z$ are adjacent, then they have no common neighbour, by the assumptions of the second and 
third bullets above; so contracting $zz'=e$ (say) will not make any parallel edges. 
Let $G'$ be the graph obtained from $G$ by contracting $e$ into a new vertex $z''$ say, and let 
$Z'=(Z\setminus \{z,z'\})\cup \{z''\}$. Then it is easy to see that
\begin{itemize}
\item $(G',Z')$ is a plantation; 
\item every $Z$-covering path of $(G,Z)$ contains $e$; so for every $Z$-covering path $P$ of $(G,Z)$, there is a $Z'$-covering 
path $P'$ of $(G',Z')$ with $E(P')=E(P)\cup \{e\}$; and
\item for every  $Z'$-covering path $P'$ of $(G',Z')$, there is at most one $Z$-covering path $P$ of $(G,Z)$ with $E(P')=E(P)\cup \{e\}$.
\end{itemize}
Consequently, in this case $n(G,Z)\le n(G',Z')$ and we can again win by induction on $|G|$. This proves \ref{simplify}.~\bbox

A {\em multiset} is a set together with a positive integer assigned to each member of the set, called its {\em multiplicity}.
The next result implies that if $(G,Z)$ is dyadic, every $Z$-covering path $P$ is determined by the set of edges of $P$ with an 
end in $Z$. A {\em linear forest} is a forest in which every component is a path; and the {\em end-multiset}
of a linear forest $H$ is the multiset of ends of the components of $H$, where an end of a component $P$ of $H$
has multiplicity one if $E(P)\ne \emptyset$, and multiplicity two if $E(P)=\emptyset$.

\begin{thm}\label{Tjoin}
Let $F$ be a forest, and let $X$ be a multiset of vertices of $V(F)$. Then there is at most one linear forest that is a subgraph of $F$ with end-multiset equal to $X$.
\end{thm}
\Proof
We proceed by induction on $|V(F)|$. If some vertex in $X$ has multiplicity at least three in $X$, then there is no linear forest
with end-multiset $X$. If  some vertex $v$ in $X$ has multiplicity two in $X$, then $v$ is a component of 
every linear forest in $F$ with end-multiset $X$, so the result follows by deleting $v$. Hence we may assume that 
every vertex in $X$ has multiplicity one. Also, from the inductive hypothesis applied to each component, we may assume 
that $F$ is connected. A {\em leaf} of $F$ means a vertex with degree one in $F$. 
If some leaf of $F$ is not in $X$, we may delete it and apply the inductive hypothesis, so we assume
all leaves of $F$ belong to $X$.
If $F$ is a path, the result is clear, so we assume $F$ is not a path.
Let us say a {\em shoot} of $F$ is a path of $F$ with one end a leaf of $F$, such that all its internal vertices have degree 
two in $F$, and maximal with both these properties. Every shoot has length at least one, one of its ends is a leaf of $F$,
and the other has degree at least three in $F$, from the maximality of the shoot and since $F$ is not a path. (Let us call the 
end of degree at least three the {\em inner end}.) Let $F'$
be obtained from $F$ by deleting all vertices of $F$ that belong to shoots and have degree at most two in $F$. Then $F'$ is non-null,
and therefore a tree; let $u$ be a vertex of $F'$ with degree at most one in $F'$. Since $u$ is not a leaf
of $F$, it is the inner end of some shoot of $F$; and therefore it has degree at least three in $F$; and so is the inner end of 
at least two shoots of $F$, say $P,P'$. But then $P\cup P'$ is a component of every linear forest in $F$ with end-multiset $X$,
and the result follows from the inductive  hypothesis by deleting $V(P\cup P')$. This proves \ref{Tjoin}.~\bbox

\bigskip

We will show:
\begin{thm}\label{finalcount}
Let 
\begin{align*}
d_1&=6( 2\phi(s)+7s-4)+4s\cdot s!\\
d_2&=8 \cdot s!(2\cdot s!+s)\phi(s)+8\\ 
d_3&=4s\cdot s!.
\end{align*}
If $(G,Z)$ is a  plantation, then $n(G,Z)\le |G|^{d_1}2^{d_2|Z|+d_3}$.
\end{thm}
\Proof By \ref{simplify} we may assume that $(G,Z)$ is dyadic.
Let $\delta_G(Z)$ be the set of edges of $G$ between $Z$ and $V(G)\setminus Z$.
\\
\\
(1) {\em For each subset $D$ of $\delta_G(Z)$, there is at most one $Z$-covering path $P$ with $E(P)\cap \delta_G(Z)=D$.}
\\
\\
To see this, let $X$ be the set of vertices in $V(G)\setminus Z$ incident with a vertex in $D$, made into a multiset by
declaring that the multiplicity of
a vertex $v$ in $X$
is the number of edges in $D$ incident with $v$. If $P$ is a $Z$-covering path with $E(P)\cap \delta_G(Z)=D$,
then $P\setminus Z$ is a linear forest with end-multiset $X$, and so $P$ is unique by \ref{Tjoin}. This proves (1).

\bigskip

Thus, in order to bound $n(G,Z)$, it is enough to bound the number of different intersections of such paths with $\delta_G(Z)$,
and we will use \ref{countedges} to do this.
Let $F=G\setminus Z$.
By \ref{countedges},
there exist $X\subseteq Z$ with $|X|\le 2\phi(s)+7s-4$  and $Y\subseteq V(G)\setminus Z$
with $|Y|\le  2s\cdot s!$ and with the following property. Let $N(X)$ be the set of vertices of $F$ with a neighbour in $X$.
For $i = 1,2$, let $N_i$ be the set of all $v\in V(F)\setminus (Y\cup  N(X))$ that have exactly $i$ neighbours in $Z$;
and let $N_0$ be the set of all $v\in N_1$ such that the component of $F\setminus (Y\cup N(X)\cup N_2)$
containing $v$ contains no other vertex in $N_1$.
There are at most
$d_2|Z|+d_3$ 
edges between $Z\setminus X$ and $V(F)\setminus (N(X)\cup N_0)$.

The edges of $\delta_G(Z)$  fall into three groups that we will handle differently, as follows:
\begin{itemize}
\item {\bf Edges between $\boldsymbol{Z}$ and $\boldsymbol{N(X)}$.} If $P$ is a $Z$-covering path, then every edge of $P$ between $Z$ and $N(X)$
belongs to a two-edge subpath of $P$ with an end in $X$. There are at most $2|X|$ such subpaths in $P$, and for each $x\in X$
the number of two-edge paths in $G$ with one end $x$ is at most $|G|^2$. Thus the number of possibilities for the set of edges
of $P$ between $Z$ and $N(X)$ is at most $|G|^{4|X|}$.
\item {\bf Edges between $\boldsymbol{Z\setminus X}$ and $\boldsymbol{N_0}$.} Let $T_1\LL T_k$ be the components of $F\setminus (Y\cup N(X)\cup N_2)$
that contain a unique vertex in $N_1$. We claim that if $P$ is a $Z$-covering path, there are at most $d_1$ values of $i\in \{1\LL k\}$
such that $P$ contains the edge between $Z$ and $V(T_i)$. To see this, suppose that $P$ contains the unique 
edge between $Z$ and $V(T_i)$.
Since both ends of $P$ are in $Z$, $P$ contains at least one edge between $V(T_i)$ and $V(F)\setminus V(T_i)$, say $uv$,
where $v\in V(F)\setminus V(T_i)$. Since $T_i$ is a component of $F\setminus (Y\cup N(X)\cup N_2)$, it follows that
$v\in Y\cup N(X)\cup N_2$. 
Suppose that $v\in N_2$; then $v\in V(P)$, but the two neighbours of $v$ in $Z$ also belong to $V(P)$,
and so $v$ has degree more than two in $P$, a contradiction. Thus $v\in Y\cup N(X)$.
We have shown then that the number of $i$ such that $P$ contains the unique 
edge between $Z$ and $V(T_i)$ is at most the number of edges of $P$ between $V(T_1\cupcup T_k)$ and $Y\cup N(X)$.
For each $v\in Y$ there are at most two edges  of $P$ between $V(T_1\cupcup T_k)$ and $v$; and for
each $v\in N(X)$ there is at most one such edge, since there is an edge of $P$ between $v$ and $X$. Since at most $2|X|$
vertices of $P$ belong to $N(X)$, it follows that 
there are at most $2|X|+2|Y|$ edges of $P$ between $V(T_1\cupcup T_k)$ and $Y\cup N(X)$.
Consequently $P$ contains at most $2|X|+2|Y|$ edges between $Z\setminus X$ and $N_0$. There are at most 
$|G|$ edges between $Z\setminus X$ and $N_0$, and so there are at most $|G|^{2|X|+2|Y|}$ possibilities for the subset that belongs to $P$.
\item {\bf Edges between $\boldsymbol{Z\setminus X}$ and $\boldsymbol{V(F)\setminus (N(X)\cup N_0)}$.} From the choice of $X,Y$, there are only $d_2|Z|+d_3$ such 
edges, so the number of possibilities for the subset that belongs to a $Z$-covering path is at most $2^{d_2|Z|+d_3}$.

\end{itemize}
It follows that the number of possibilities for $E(P)\cap \delta_G(Z)$ is at most the product of these three; and so
$$n(G,Z)\le |G|^{4|X|}|G|^{2|X|+2|Y|}2^{d_2|Z|+d_3}\le |G|^{d_1} 2^{d_2|Z|+d_3}.$$
This proves \ref{finalcount}.~\bbox

We deduce \ref{mainthm}, which we restate:
\begin{thm}\label{count}
For all integers $s\ge 1$, there exist $c_1,c_2,c_3\ge 0$ such that 
if $G$ is $s\mathcal{O}$-free, and $Z\subseteq V(G)$ is
a cycle-hitting set, then $G$ has at most
$|G|^{c_1}2^{c_2|Z|+c_3}$ induced paths.
\end{thm}
\Proof
There are at most $|G|^2$ induced paths that are vertex-disjoint from $Z$, since such paths are determined by their ends.
Let us count the induced paths that have a vertex in $Z$. For each such path $Q$, with ends $s,t$ say,
 let $a$ be the vertex of $Q$ in $Z$ that is closest to $s$ in $Q$, and define $b$ similarly for $t$. (Possibly $s=a$, or $a=b$, or $b=t$.)
Thus $Q$ is divided
into three subpaths: the subpath between $s$ and $a$, the subpath between $a$ and $b$, and the subpath between $b$ and $t$.
There are only $|G|^2/2$ possibilities for the first part, since it is determined by its first vertex and penultimate vertex;
and similarly there are only $|G|^2/2$ possibilities for the last part. We need to count the possibilities for the 
middle part $P$ say, between $a$ and $b$. Let $Z'=Z\cap V(P)$; then $P$ is a $Z'$-covering path in the plantation
$(G\setminus (Z\setminus Z'), Z')$, and so, with $d_1,d_2,d_3$ as in \ref{finalcount}, for each choice of $Z'$,
the number of choices of $P$ is at most $|G|^{d_1}2^{d_2|Z|+d_3}$.
Since there are only $2^{|Z|}-1$ choices for $Z'$ (since $Z'\ne \emptyset$), there are only $|G|^{d_1}2^{d_2|Z|+d_3}(2^{|Z|}-1)$ choices for $P$ in total, and hence
only 
$$|G|^{d_1+4}2^{d_2|Z|+d_3}(2^{|Z|}-1)+|G|^2\le |G|^{d_1+4}2^{(d_2+1)|Z|+d_3}$$ 
choices for $Q$. 
This proves \ref{count}.~\bbox

\section*{Acknowledgement}
Thanks to Maria Chudnovsky and Sophie Spirkl for stimulating discussions, and for their work on parts of this paper.

\end{document}